\documentclass[11pt,reqno]{amsart}
\usepackage[left=2.5cm,right=2.5cm,top=2.5cm, bottom=2.5cm]{geometry}
\usepackage{amsfonts,epsfig,fancyhdr}
\usepackage{amsmath,amssymb,amsthm, mathrsfs, fancyhdr,latexsym,graphicx,graphics}

\usepackage{tikz-cd}

\usepackage{pdfsync}
\usepackage{url,hyperref}
\usepackage{euscript}
\usepackage{wrapfig}
\usepackage{caption}
\usepackage{tcolorbox}
\usepackage{graphicx}
\numberwithin{equation}{section}

\subjclass[2010]{Primary ...; ...}

\numberwithin{equation}{section}

\newtheorem{theorem}{Theorem}[section]
\newtheorem{cor}[theorem]{Corollary}

\newtheorem{rem}[theorem]{Remark}

\newtheorem{ejem}[theorem]{Example}

\newtheorem{prop}[theorem]{Proposition}
\newtheorem{thm}[theorem]{Theorem}

\theoremstyle{definition}
\newtheorem{definition}{Definition}[section]

\usepackage[dvipsnames]{xcolor}
\usepackage{enumerate}

\title[Geometric PDEs-HT for Maps into F$_\text{man}$-Lie Groups]{Geometric Partial Differential Equations of Hydrodynamic Type for Maps into F$_\text{man}$-Lie Groups}

\author[S. Castañeda-Montoya]{Santiago Castañeda-Montoya}
\address{Instituto de Matemáticas, FCEyN, Universidad de Antioquia, 50010  Medellín, Colombia}
\email{santiago.castanedam@udea.edu.co}

\author[A. Torres-Gomez]{Alexander Torres-Gomez}
\email{galexander.torres@udea.edu.co}


\begin{document}

\maketitle

\begin{abstract}
We develop a geometric framework for partial differential equations of hydrodynamic type on maps from a two-dimensional manifold into an F$_\text{man}$-Lie group. Using the left-invariant structure and the Maurer–Cartan form, we reduce these equations to non-linear evolution systems on the underlying F$_\text{man}$-algebra, governed by its Lie bracket and commutative associative product. We investigate the role of Poisson subalgebras, establishing that flatness of the canonical connection forces two-step nilpotency and yields a linearized evolution equation upon quotienting by the center. Furthermore, we introduce the notion of a Poisson splitting and derive the corresponding tangential and transversal evolution equations. Under suitable reductivity assumptions, we show that their non-linear coupling is governed by the curvature of the canonical principal connection on the homogeneous space $G/H$. These results reveal the interplay between the algebraic structure of F$_\text{man}$-Lie groups and the geometric dynamics of hydrodynamic-type equations.
\end{abstract}

\tableofcontents

%%%%%%%%%%%%%%%%%%%%%%%%%%%%%%%%%%%%%%%%% Section
\section{Introduction}
In \cite{dubrovin-novikov}, Dubrovin and Novikov introduced a general theory of \emph{systems of hydrodynamic type}, a class of first-order quasilinear partial differential equations. Their work established a profound synthesis connecting flat pseudo-Riemannian geometry, Hamiltonian mechanics (via Poisson brackets of hydrodynamic type), Riemann invariants, and integrability. Furthermore, this foundational framework anticipated the later geometric theory of Frobenius manifolds \cite{dubrovin}.\\

Through the study of Frobenius manifolds \cite{dubrovin}, Dubrovin provided a crucial geometric perspective on the Witten–Dijkgraaf–Verlinde–Verlinde (WDVV) equations arising in two-dimensional topological field theories. Generalizing this geometric framework, Hertling and Manin \cite{hertling-manin} introduced the concept of a weak Frobenius manifold, nowadays commonly referred to as an \emph{F-manifold}. This construction relaxes several conditions inherent to Frobenius manifolds. Notably, every Frobenius manifold is, in particular, an F-manifold.\\

Roughly speaking, an F-manifold consists of a smooth manifold $M$ whose $C^\infty(M)$-module of vector fields, $\mathfrak{X}(M)$, is endowed with a commutative and associative product $\circ$. This product is subject to an integrability constraint known as the Hertling–Manin condition, which links the Lie bracket and the multiplication $\circ$ through a specific cubic identity (see Definition \ref{fmanifold}). F-manifolds are ubiquitous, appearing in diverse areas such as singularity theory, quantum cohomology, symplectic geometry, and integrable systems \cite{manin1999, hertling2002}. A setting of primary interest for the present work is established by Lorenzoni, Pedroni, and Raimondo \cite{LPR}, who investigated the interplay between F-manifolds with compatible connections and integrable systems of hydrodynamic type.\\

An important algebraic simplification occurs when the underlying manifold $M$ of an F-manifold is a Lie group $G$, a structure we termed an \emph{F$_{\text{man}}$-Lie group} \cite{our-paper}. We develop a global, coordinate-free geometric formulation of systems of partial differential equations of hydrodynamic type by realizing their solutions as smooth maps—termed hydrodynamic maps—from a two-dimensional manifold into an F$_\text{man}$-Lie group $G$. Using left trivialization via the Maurer–Cartan form, we show that a hydrodynamic map is locally governed by an evolutionary equation on the associated F$_\text{man}$-algebra $(\mathfrak{g}, [\cdot, \cdot], \circ)$ of the form $P_t - a \circ P_x + [a \circ P, P] = 0$, unifying the commutative F-manifold multiplication with the zero-curvature geometry of the underlying Lie group. We characterize compatible and flat left-invariant connections on $G$ in purely algebraic terms and prove that the canonical F$_\text{man}$-connection is compatible with the multiplication tensor if and only if the algebraic Leibnizator is totally symmetric. Furthermore, we show that flat Poisson subalgebras induce 2-step nilpotent foliations along which the hydrodynamic flow linearizes modulo the center. Finally, we introduce the notion of a Poisson splitting $\mathfrak{g} = \mathfrak{h} \oplus \mathfrak{m}$ and prove both local and global bundle splitting theorems, demonstrating that the dynamic coupling between the vertical flow along Poisson leaves and the horizontal transport across the homogeneous base $G/H$ is governed by the curvature $2$-form of the canonical principal connection.\\

The remainder of this paper is structured as follows. Section \ref{prelim} provides a brief review of PDEs of hydrodynamic type, F-manifolds, F$_{\text{man}}$-Lie groups, and their associated F$_{\text{man}}$-algebras. Section \ref{geometric-pde} is dedicated to constructing the coordinate-free geometrization of hydrodynamic PDEs, treating the independent variables as local coordinates on a two-dimensional manifold and the dependent field functions as coordinates on an F$_{\text{man}}$-Lie group. Finally, Section \ref{poisson-decomp} develops a structural splitting of these hydrodynamic PDEs into a Poisson subalgebra and a transverse vector subspace.

%%%%%%%%%%%%%%%%%%%%%%%%%%%%%%%%%%%%%%%%% Section

\section{Preliminaries}\label{prelim}

In this section, we review the fundamental algebraic and geometric concepts that underpin our work. We begin with a brief recap of partial differential equations (PDEs) of hydrodynamic type. Next, we formulate F-manifolds using a symmetric tensor and the Leibnizator, specialize this framework to F$_{\text{man}}$-Lie groups and their associated F$_{\text{man}}$-algebras, and review compatible connections as well as Poisson algebra distributions.

\subsection{PDEs of Hydrodynamic Type} 
A system of $n$ first-order quasilinear partial differential equations in one spatial variable $x$ and one temporal variable $t$ is said to be of \emph{hydrodynamic type} if it takes the form
\begin{equation}\label{hydro-eq-coord-full}
\frac{\partial u^i(x,t)}{\partial t} - V^i{}_j(u) \frac{\partial u^j(x,t)}{\partial x} = 0, \qquad i = 1, \dots, n,
\end{equation}
where $(x,t)$ are local coordinates on an open subset $D \subseteq \mathbb{R}^2$, the map $u = (u^1, \dots, u^n) \colon D \to \mathbb{R}^n$ is smooth, and $V^i{}_j(u)$ is a matrix-valued function that depends explicitly on $u$, but not on its partial derivatives. For brevity, we write \eqref{hydro-eq-coord-full} in the compact notation
\begin{equation}\label{hydro-eq-coord}
u_t^i - V^i{}_j(u) , u_x^j = 0.
\end{equation}
The term \emph{hydrodynamic type} originates from the fact that equations of this form arise naturally as dispersionless (semi-classical) limits of physical field theories and fluid dynamics. In what follows, we refer to such a system as a \emph{PDE-HT}. Below, we present a coordinate-free, geometric formulation of these systems.

\subsection{F-manifolds}
The geometric structure of F-manifolds \cite{hertling-manin, hertling2002, manin1999} relies on endowing the space of vector fields with a commutative and associative product, which can be cleanly formulated using tensor fields.

\begin{definition}[Commutative Product and Tensor $S$]
Let $M$ be a smooth manifold and $\mathfrak{X}(M)$ the $\mathcal{C}^\infty(M)$-module of smooth vector fields. A smooth, commutative and associative multiplication $\circ \colon \mathfrak{X}(M) \times \mathfrak{X}(M) \to \mathfrak{X}(M)$ admits a geometric representation as a symmetric $(1,2)$-tensor field $S \in \Gamma(T^{(1,2)}M)$ defined by:
\begin{equation*}
S(Y, Z) := Y \circ Z, \quad \text{for all } Y, Z \in \mathfrak{X}(M).
\end{equation*}
The symmetry of $S$ follows directly from the commutativity of $\circ$.
\end{definition}

To measure the obstruction to the Lie bracket acting as a derivation on the commutative and associative $\mathbb{R}$-algebra $(\mathfrak{X}(M), \circ)$, one introduces the \emph{Leibnizator}.
\begin{definition}[Leibnizator]
For any triplet of vector fields $X, Y, Z \in \mathfrak{X}(M)$, the \emph{Leibnizator} is the trilinear map $L \colon \mathfrak{X}(M) \times \mathfrak{X}(M) \times \mathfrak{X}(M) \to \mathfrak{X}(M)$ given by:
\begin{equation*}
L(X, Y, Z) := [X, Y \circ Z] - [X, Y] \circ Z - Y \circ [X, Z].
\end{equation*}
\end{definition}

\begin{rem}
The Leibnizator $L$ is $\mathcal{C}^\infty(M)$-multilinear with respect to its second and third arguments, while acting as a first-order differential operator in its first argument due to the properties of the Lie bracket. Furthermore, $L$ corresponds precisely to the Lie derivative of the symmetric tensor $S$ along the vector field $X$ \cite[Lemma 2.5]{our-paper}:
\begin{equation*}
(\mathcal{L}_X S)(Y, Z) = L(X, Y, Z),
\end{equation*}
where $\mathcal{L}$ denotes the Lie derivative of tensor fields.
\end{rem}\begin{definition}[F-Manifold]\label{fmanifold}
An \emph{F-manifold} is a pair $(M, \circ)$ (or equivalently $(M, S)$) consisting of a smooth manifold $M$ and a commutative and associative product $\circ$ satisfying the \emph{Hertling–Manin condition}:
\begin{equation}\label{h-m-condition}
L(X \circ Y, Z, W) = L(X, Z, W) \circ Y + X \circ L(Y, Z, W), \quad \text{for all } X, Y, Z, W \in \mathfrak{X}(M).
\end{equation}
\end{definition}

\subsection{F$_\text{man}$-Lie Groups and F$_\text{man}$-Algebras}
The theory of F-manifolds becomes more tractable when the underlying manifold is a Lie group $G$ and all geometric structures are invariant under left translations. For a detailed treatment of F$_{\text{man}}$-Lie groups and their corresponding F$_{\text{man}}$-algebras, we refer the reader to \cite{our-paper} (see also \cite{TGV}).

\begin{definition}[F$_{\text{man}}$-Lie Group]
An \emph{F$_{\text{man}}$-Lie group} $(G, \circ)$ is a Lie group $G$ equipped with an F-manifold structure such that the multiplication $\circ$ (or equivalently, the symmetric tensor field $S$) on the tangent bundle $TG$ is left-invariant.
\end{definition}

Recall that there is a canonical vector space isomorphism between the space of left-invariant vector fields on $G$, denoted by $\mathfrak{X}^+(G)$, and the tangent space at the identity element $e \in G$, given by the Lie algebra $\mathfrak{g} \cong T_e G$. We identify each left-invariant vector field $a^+$ with its value at the identity, $a := a^+\vert{}_e \in \mathfrak{g}$. Since the Lie bracket is inherently left-invariant and the product $\circ$ is left-invariant by assumption, the Leibnizator $L$ constitutes a left-invariant trilinear map on $\mathfrak{X}^+(G)$. Consequently, $L$ is uniquely determined by its evaluation at $e \in G$. Restricting this operation to the Lie algebra $\mathfrak{g} \cong T_e G$ yields the \emph{algebraic Leibnizator}:
\begin{equation*}
L(a, b, c) = [a, b \circ c] - [a, b] \circ c - b \circ [a, c], \quad \text{for all } a, b, c \in \mathfrak{g}.
\end{equation*}
In this algebraic context, the Lie bracket $[\cdot, \cdot]$ and the product $\circ$ on $\mathfrak{g}$ are the operations induced by evaluating the corresponding left-invariant vector fields at the identity:
\begin{equation*}
[a, b] := [a^+, b^+]|_e \quad \text{and} \quad a \circ b := (a^+ \circ b^+)|_e.
\end{equation*}

\begin{definition}[F$_{\text{man}}$-Algebra]
An \emph{F$_{\text{man}}$-algebra} is a finite-dimensional real vector space $\mathfrak{g}$ equipped with a Lie bracket $[\cdot, \cdot]$ and a commutative and associative product $\circ$ satisfying the \emph{algebraic Hertling–Manin condition} via its algebraic Leibnizator:
\begin{equation*}
L(a \circ b, c, d) = L(a, c, d) \circ b + a \circ L(b, c, d), \quad \text{for all } a, b, c, d \in \mathfrak{g}.
\end{equation*}
\end{definition}

Every F$_{\text{man}}$-Lie group naturally induces an F$_{\text{man}}$-algebra structure on its Lie algebra $\mathfrak{g}$.

\subsubsection{Poisson Algebra Distributions and Foliations}
Within an F$_{\text{man}}$-algebra, certain vector subspaces display a deeper compatibility between the Lie bracket and the commutative product.

\begin{rem}[Poisson Subalgebra]
Let $(\mathfrak{g}, [\cdot, \cdot], \circ)$ be an F$_{\text{man}}$-algebra. A subspace $\mathfrak{h} \subseteq \mathfrak{g}$ that is closed under both the Lie bracket $[\cdot, \cdot]$ and the product $\circ$, and on which the Leibnizator vanishes identically (i.e., $L(a, b, c) = 0$ for all $a, b, c \in \mathfrak{h}$), is called a \emph{Poisson subalgebra}. As is well known, on a Poisson algebra $(\mathfrak{h}, [\cdot, \cdot], \circ)$, the Lie bracket acts as a derivation over the commutative and associative algebra $(\mathfrak{h}, \circ)$.
\end{rem}

Geometrically, this algebraic property dictates the integrability of the corresponding left-invariant distributions on the Lie group $G$ \cite{our-paper}.

\begin{theorem}[Integrability and Foliation]
Let $G$ be an F$_{\text{man}}$-Lie group and $\mathfrak{h} \subseteq \mathfrak{g}$ a Poisson subalgebra. The left-invariant distribution $\mathcal{D}_{\text{Pois}}$ generated by $\mathfrak{h}$ is completely integrable. Consequently, $\mathcal{D}_{\text{Pois}}$ defines a regular foliation of $G$, where each leaf locally inherits a Poisson structure.
\end{theorem}
\begin{proof}
See \cite[Theorem 3.13]{our-paper}.
\end{proof}

%%%%%%%%%%%%%%%%%%%%%%%%%%%%%%%%%%%%%%% Section

\section{Geometric PDE of Hydrodynamic Type}\label{geometric-pde}
In this section, we present a coordinate-free, geometric formulation of partial differential equations of hydrodynamic type (PDEs-HT). Rather than treating the dependent variables merely as local scalar functions, we formulate a prototype geometric PDE-HT in terms of a smooth map $\phi \colon \Sigma \to M$ from a two-dimensional manifold $\Sigma$ (an abstract surface representing the space-time domain) into an $n$-dimensional target manifold $M$. We refer to such a map $\phi$ as a \emph{hydrodynamic map}.

\subsection{Coordinate-free Formulation}
Let $\Sigma$ be a smooth two-dimensional manifold (an abstract surface) and $M$ an $n$-dimensional target manifold. Consider a smooth map $\phi \colon \Sigma \to M$. The differential (or pushforward) of $\phi$ is a smooth bundle map between the corresponding tangent bundles:
\begin{equation*}
d\phi \colon T\Sigma \longrightarrow TM.
\end{equation*}
For any point $p \in \Sigma$, the restriction of $d\phi$ to the tangent space $T_p \Sigma$ yields a linear map between vector spaces:
\begin{equation*}
d\phi_p \colon T_p \Sigma \longrightarrow T_{\phi(p)} M.
\end{equation*}

Choosing local coordinates $(x,t)$ on $\Sigma$ induces a canonical local vector field frame for $T\Sigma$, given by $\frac{\partial}{\partial x}, \frac{\partial}{\partial t} \in \Gamma(T\Sigma)$. Evaluating $d\phi$ on this coordinate frame yields two geometric assignments at each point $p \in \Sigma$:
\begin{equation*}
\phi_x(p) := d\phi_p \left( \left. \frac{\partial}{\partial x} \right|_p \right) \in T_{\phi(p)} M, \qquad
\phi_t(p) := d\phi_p \left( \left. \frac{\partial}{\partial t} \right|_p \right) \in T_{\phi(p)} M.
\end{equation*}

To treat $\phi_x$ and $\phi_t$ rigorously as standard sections over the domain manifold $\Sigma$, we construct the pullback bundle $\phi^* TM \to \Sigma$. The pullback bundle $\phi^* TM$ is the fiber bundle over $\Sigma$ whose fiber at $p \in \Sigma$ is defined by $(\phi^* TM)_p := T_{\phi(p)} M$. As a total space, it is defined via the fiber product
\begin{equation*}
\phi^* TM = \{ (p, v) \in \Sigma \times TM \mid \phi(p) = \pi_M(v) \},
\end{equation*}
equipped with the canonical projection $\pi_{\phi} \colon \phi^* TM \to \Sigma$ given by $\pi_{\phi}(p, v) = p$. More explicitly, this construction fits into the commutative diagram:
\[
\begin{tikzcd}
    \phi^*TM \arrow[r] \arrow[d, "\pi_\phi"'] & TM \arrow[d, "\pi_M"] \\
    \Sigma \arrow[r, "\phi"'] & M
\end{tikzcd}
\]

Since $d\phi_p \colon T_p \Sigma \to (\phi^* TM)_p = T_{\phi(p)} M$ at every point $p \in \Sigma$, the differential $d\phi$ is naturally a smooth section of the bundle $T^*\Sigma \otimes \phi^* TM$. That is,
\begin{equation*}
d\phi \in \Gamma(T^* \Sigma \otimes \phi^*TM) = \Omega^1(\Sigma, \phi^* TM).
\end{equation*}
Evaluating $d\phi \in \Omega^1(\Sigma, \phi^* TM)$ on the coordinate vector fields $\partial_x, \partial_t \in \Gamma(T\Sigma)$ yields the smooth sections $\phi_x, \phi_t \in \Gamma(\phi^* TM)$.\\

Within this bundle framework, a system of hydrodynamic type admits a coordinate-free, intrinsic formulation. Let $V \in \Gamma(\operatorname{End}(TM))$ be an endomorphism field on $TM$, so that at each point $q \in M$, $V_q \colon T_q M \to T_q M$ is a linear endomorphism. Pulling $V$ back along $\phi$ yields a smooth bundle endomorphism $\phi^* V \colon \phi^* TM \to \phi^* TM$, which at each point $p \in \Sigma$ acts as
\begin{equation*}
(\phi^* V)_p = V_{\phi(p)} \colon T_{\phi(p)} M \longrightarrow T_{\phi(p)} M.
\end{equation*}

Consequently, the composition $(\phi^* V)(\phi_x)$ is well-defined: both $\phi_t$ and $(\phi^* V)(\phi_x)$ are sections of the exact same vector bundle $\phi^* TM$ over $\Sigma$. The geometric PDE of hydrodynamic type is then formulated as the vanishing section condition
\begin{equation}\label{eq:hydro-bundle}
\phi_t - (\phi^* V)(\phi_x) = 0 \quad \in \Gamma(\phi^* TM),
\end{equation}
or equivalently,
\begin{equation*}
d\phi(\partial_t) - (\phi^* V)\bigl( d\phi(\partial_x) \bigr) = 0.
\end{equation*}

In local coordinates $(u^1, \dots, u^n)$ on $M$, equation \eqref{eq:hydro-bundle} expands as
\begin{equation*}
\sum_{i=1}^n \left( u_t^i - V^i{}_j(u) \, u_x^j \right) \left. \frac{\partial}{\partial u^i} \right|_{\phi} = 0.
\end{equation*}
Since the local coordinate frame fields $\left. \frac{\partial}{\partial u^i} \right\vert{}_{\phi}$ are pointwise linearly independent, this reduces directly to
\begin{equation*}
u_t^i - V^i{}_j(u) \, u_x^j = 0, \qquad i = 1, \dots, n,
\end{equation*}
recovering the standard coordinate system of hydrodynamic type.

\subsection{F$_\text{man}$-Lie group as a Target Manifold}
Let $(G, \circ)$ be an F$_{\text{man}}$-Lie group and $\phi \colon \Sigma \to G$ a smooth map from a smooth surface $\Sigma$ into $G$. Choosing local coordinates $(x,t)$ on $\Sigma$, we obtain the pushed-forward coordinate frame sections $\phi_x, \phi_t \in \Gamma(\phi^*TG)$. The hydrodynamic partial differential equation associated with a vector field $X \in \mathfrak{X}(G)$ is given by:
\begin{equation}\label{HydroPDE-G-VF}
\phi_t = X(\phi) \circ \phi_x.
\end{equation}
This framework is particularly well-suited to our setting when $X$ is chosen to be a left-invariant vector field on $G$. The Lie group structure of $G$ allows for an even more intrinsic formulation by left-trivializing all geometric quantities via the Maurer–Cartan form.

\subsubsection{Maurer-Cartan Form and Trivialization}
We begin by recalling two fundamental integration theorems for Lie-algebra-valued differential forms (see Sharpe \cite[Chapter 3, \S 6, Theorems 6.1 \& 7.14]{Sharpe}).

\begin{thm}[Local Version]
Let $G$ be a Lie group with Lie algebra $\mathfrak{g}$. Let $A$ be a $\mathfrak{g}$-valued $1$-form on a smooth manifold $\Sigma$ satisfying the Maurer–Cartan zero-curvature equation:
\begin{equation*}
dA + \tfrac{1}{2}[A, A] = 0.
\end{equation*}
Then, for each point $p \in \Sigma$, there exist a neighborhood $U \subseteq \Sigma$ of $p$ and a smooth map $\phi \colon U \to G$ such that $A\vert{}_U = \phi^* \theta$, where $\theta$ denotes the left Maurer–Cartan form on $G$.
\end{thm}\begin{thm}[Global Version]
Let $G$ be a Lie group with Lie algebra $\mathfrak{g}$, let $\Sigma$ be a smooth connected manifold, and let $A \in \Omega^1(\Sigma, \mathfrak{g})$. Then $A$ is the pullback $A = \phi^*\theta$ of some smooth map $\phi \colon \Sigma \to G$ if and only if:
\begin{enumerate}[(i)]
\item $dA + \tfrac{1}{2}[A, A] = 0$, and
\item the monodromy representation $\Phi_A \colon \pi_1(\Sigma, p) \to G$ is trivial.
\end{enumerate}
Moreover, when these conditions hold, the generating map $\phi \colon \Sigma \to G$ is unique up to left translation by a constant element of $G$.
\end{thm}

The Lie group geometry allows us to canonically identify each tangent space $T_g G$ with the Lie algebra $\mathfrak{g} \cong T_e G$ via left translation. For any tangent vector $v_g \in T_g G$, its left-trivialization is defined by $(dL_{g^{-1}})_g(v_g) \in \mathfrak{g}$. Consequently, for a smooth map $\phi \colon \Sigma \to G$, we define the components:
\begin{equation*}
A_t := (dL_{\phi^{-1}})_\phi(\phi_t) \in \mathfrak{g}, \qquad A_x := (dL_{\phi^{-1}})_\phi(\phi_x) \in \mathfrak{g}.
\end{equation*}
These components define smooth $\mathfrak{g}$-valued functions $A_t, A_x \colon \Sigma \to \mathfrak{g}$, or equivalently, the pulled-back Maurer–Cartan $1$-form:
\begin{equation*}
A = \phi^{-1} d\phi = \phi^* \theta \in \Omega^1(\Sigma, \mathfrak{g}).
\end{equation*}
In local coordinates $(x,t)$ on $\Sigma$, this $\mathfrak{g}$-valued $1$-form expands as:
\begin{equation*}
A = A_t \, dt + A_x \, dx.
\end{equation*}

Now, let $a^+$ be a left-invariant vector field on $G$, and let $a := a^+\vert{}_e \in \mathfrak{g}$ be its corresponding Lie algebra element. Because the product $\circ$ is left-invariant, it restricts to a commutative and associative bilinear multiplication $\circ \colon \mathfrak{g} \times \mathfrak{g} \to \mathfrak{g}$ on the Lie algebra. Left-trivializing equation \eqref{HydroPDE-G-VF} yields:
\begin{equation*}
(dL_{\phi^{-1}})_\phi \left( a^+(\phi) \circ \phi_x \right) = a \circ A_x.
\end{equation*}
Thus, the hydrodynamic relation simplifies dramatically to the linear algebraic condition:
\begin{equation*}
A_t = a \circ A_x.
\end{equation*}
In other words, the nonlinear-looking geometric PDE \eqref{HydroPDE-G-VF} with target $G$ reduces under left-trivialization to:
\begin{equation*}
A_t = S_a(A_x),
\end{equation*}
where $S_a \in \operatorname{End}(\mathfrak{g})$ is the endomorphism given by left multiplication $S_a(b) := a \circ b$.\\

However, not every $\mathfrak{g}$-valued $1$-form $A$ originates from a map $\phi \colon \Sigma \to G$. Because $A = \phi^* \theta$, it must satisfy the integrability condition given by the Maurer–Cartan structural equation:
\begin{equation*}
dA + \frac{1}{2}[A, A] = 0,
\end{equation*}
which in local coordinates $(x,t)$ takes the form:
\begin{equation*}
\partial_t A_x - \partial_x A_t + [A_t, A_x] = 0.
\end{equation*}
Thus, the local hydrodynamic-map system into an F$_{\text{man}}$-Lie group translates into the fundamental coupled system:
$$ \begin{cases} A_t=a\circ A_x,\\[4pt] \partial_tA_x-\partial_xA_t+[A_t,A_x]=0. \end{cases}$$
Here, the first equation reflects the algebraic F-manifold product, while the second reflects the differential geometry of the underlying Lie group.\\

Substituting $A_t = a \circ A_x$ into the Maurer–Cartan equation yields:
\begin{equation*}
\partial_t A_x - \partial_x(a \circ A_x) + [a \circ A_x, A_x] = 0.
\end{equation*}
Since $a \in \mathfrak{g}$ originates from a left-invariant vector field $a^+$, it is constant over $\Sigma$. Consequently, $\partial_x(a \circ A_x) = a \circ \partial_x A_x$, giving:
\begin{equation*}
\partial_t A_x - a \circ \partial_x A_x + [a \circ A_x, A_x] = 0.
\end{equation*}
To simplify notation, set:
\begin{equation*}
P := A_x \in C^\infty(\Sigma, \mathfrak{g}), \qquad Q := A_t \in C^\infty(\Sigma, \mathfrak{g}).
\end{equation*}
The hydrodynamic system then takes the form $Q = a \circ P$ subject to the curvature constraint $P_t - Q_x + [Q, P] = 0$. Substituting $Q = a \circ P$ yields:
\begin{equation}\label{hydro.eq}
P_t - (a \circ P)_x + [a \circ P, P] = 0.
\end{equation}
For a constant element $a \in \mathfrak{g}$, this reduces to the nonlinear evolution equation:
\begin{equation}\label{hydro.eq.cons.a}
P_t - a \circ P_x + [a \circ P, P] = 0.
\end{equation}
Equation \eqref{hydro.eq.cons.a} is an evolutionary system governing the $\mathfrak{g}$-valued function $P(x,t)$, wherein the F$_{\text{man}}$-algebra structure $(\mathfrak{g}, [\cdot,\cdot], \circ)$ plays a direct and central role.

%%%%%%%%%%%%%%%%%%%%%%%%%%%%%%%%%%%%%%%%%%% Section

\section{Canonical Connection and Poisson Decomposition}\label{poisson-decomp}
We begin by recalling the general theory of left-invariant affine connections on Lie groups, which enables the reduction of geometric differential operators to purely algebraic operations on the Lie algebra. Let $G$ be a Lie group and $\mathfrak{g} = T_e G$ its associated Lie algebra, identified with the space of left-invariant vector fields on $G$.

\begin{definition}[Left-Invariant Connection]
An affine connection $\nabla$ on a Lie group $G$ is said to be \emph{left-invariant} if the left translation maps $L_g \colon G \to G$ are affine transformations for all $g \in G$. Equivalently, $\nabla$ is left-invariant if and only if for any pair of left-invariant vector fields $a^+, b^+ \in \mathfrak{X}^+(G)$, their covariant derivative $\nabla_{a^+} b^+$ is also a left-invariant vector field.
\end{definition}Due to left translation invariance, the full structure of the connection is uniquely determined by its evaluation at the identity element $e \in G$. Consequently, there is a one-to-one correspondence between left-invariant connections on $G$ and bilinear maps $\alpha \colon \mathfrak{g} \times \mathfrak{g} \to \mathfrak{g}$ on its Lie algebra (see \cite[Chapter 6, \S 3]{postnikov}), defined by:
\begin{equation*}
\alpha(a, b) = (\nabla_{a^+} b^+)|_e, \qquad \text{for all } a, b \in \mathfrak{g},
\end{equation*}
where $a^+$ and $b^+$ denote the unique left-invariant vector fields extending $a, b \in \mathfrak{g}$, i.e., $a^+\vert{}_e = a$ and $b^+\vert{}_e = b$.\\

This algebraic reduction allows fundamental geometric quantities—specifically, the torsion tensor $T$ and the Riemann curvature tensor $R$—to be expressed entirely in terms of the bilinear map $\alpha$ and the Lie bracket $[\cdot, \cdot]$ on $\mathfrak{g}$. For any elements $a, b, c \in \mathfrak{g}$, the algebraic torsion and curvature tensor evaluations $T(a, b) := T(a^+, b^+)\vert{}_e$ and $R(a, b)c := R(a^+, b^+)c^+\vert{}_e$ are given by:
\begin{align}
T(a, b) &= \alpha(a, b) - \alpha(b, a) - [a, b], \label{eq:gen_torsion} \\
R(a, b)c &= \alpha(a, \alpha(b, c)) - \alpha(b, \alpha(a, c)) - \alpha([a, b], c). \label{eq:gen_curvature}
\end{align}

\begin{rem}
In the specific setting of F$_{\text{man}}$-algebras, the canonical left-invariant connection we consider corresponds to selecting the bilinear map
\begin{equation*}
\alpha(a, b) = \frac{1}{2}[a, b] + a \circ b,
\end{equation*}
see \cite{our-paper}.
\end{rem}

\subsection{Compatible Connection}
In \cite{LPR}, Lorenzoni, Pedroni and Raimondo adapted Dubrovin’s construction of the principal hierarchy \cite{dubrovin} to F-manifolds equipped with a \emph{compatible flat connection}, a concept introduced by Manin in \cite{manin2005}. Compatible connections play a central role on F-manifolds—and specifically on F$_{\text{man}}$-Lie groups—because they provide the geometric infrastructure underlying integrable systems of hydrodynamic type. We begin by recalling the definition of a compatible connection on an F-manifold (see \cite{LPR, manin2005}).
\begin{definition}[Compatible Connection]\label{def:compatible-connection}
Let $M$ be a smooth manifold and $S \in \Gamma(T^{(1,2)}M)$ a symmetric tensor defined by a commutative vector field product $\circ$, i.e., $S(X, Y) = X \circ Y$ for all $X, Y \in \mathfrak{X}(M)$. A torsion-free affine connection $\nabla$ on $(M, S)$ is said to be \emph{compatible with $S$} (or with the product $\circ$) if the covariant derivative tensor $\nabla S$ is totally symmetric. Equivalently, for all vector fields $X, Y, Z \in \mathfrak{X}(M)$, the connection satisfies:
\begin{equation*}
(\nabla_X S)(Y, Z) = (\nabla_Y S)(X, Z).
\end{equation*}
If $\nabla$ is also flat (its Riemann curvature tensor vanishes identically), the triple $(M, S, \nabla)$ provides the geometric framework for constructing integrable hierarchies of hydrodynamic type \cite{LPR}.
\end{definition}

\begin{rem}
As Lorenzoni et al. noted in \cite[Remark 8]{LPR}, the compatibility of a flat connection $\nabla$ is equivalent to the local existence of a vector field $C$ such that, for any pair of flat vector fields $X$ and $Y$, the product $\circ$ satisfies
\begin{equation*}
X \circ Y = [X, [Y, C]].
\end{equation*}
This condition was Manin's original definition of a compatible connection \cite[Definition 2.2]{manin2005}, where the vector field $C$ was termed a \emph{local vector potential} for $\circ$. Furthermore, Hertling proved \cite[Theorem 5.2]{hertling1999} (see also \cite[Proposition 2.4]{manin2005}) that the Hertling–Manin condition \eqref{h-m-condition} is automatically satisfied whenever a manifold admits a torsion-free (not necessarily flat) connection $\nabla$ compatible with $S$. Thus, any triple $(M, S, \nabla)$ with a compatible torsion-free connection $\nabla$ is inherently an F-manifold.
\end{rem}

On an F$_{\text{man}}$-Lie group $G$, the requirements for a left-invariant connection to be compatible reduce to three algebraic conditions on the corresponding bilinear map $\alpha \colon \mathfrak{g} \times \mathfrak{g} \to \mathfrak{g}$: $\alpha$ must be torsion-free, compatible with the multiplication tensor $S$, and flat. Below, we express these requirements as explicit algebraic conditions on the F$_{\text{man}}$-algebra $(\mathfrak{g}, [\cdot, \cdot], \circ)$:

\begin{itemize}
\item \textbf{The Torsion-Free Condition.}
Requiring the torsion tensor to vanish identically ($T(a, b) = 0$) in \eqref{eq:gen_torsion} imposes a skew-symmetric constraint on the bilinear map $\alpha$. A left-invariant connection $\nabla$ on a Lie group $G$ is torsion-free if and only if its associated bilinear map $\alpha$ satisfies:
\begin{equation}\label{eq:alpha_torsion_free}
\alpha(a, b) - \alpha(b, a) = [a, b], \qquad \text{for all } a, b \in \mathfrak{g}.
\end{equation}
This shows that the skew-symmetric part of $\alpha$ is completely determined by the Lie bracket of the F$_{\text{man}}$-algebra.

\item \textbf{Compatibility with the Tensor $S$.}
Let $S(a, b) := (S(a^+, b^+))\vert{}_e = a \circ b$ for $a, b \in \mathfrak{g}$. For left-invariant vector fields $a^+, b^+, c^+ \in \mathfrak{X}^+(G)$, evaluating the covariant derivative of $S$ at the identity element yields:
\begin{equation*}
(\nabla_{a^+} S)(b^+, c^+)|_e = \alpha(a, b \circ c) - \alpha(a, b) \circ c - b \circ \alpha(a, c).
\end{equation*}
Compatibility requires $(\nabla_{a^+} S)(b^+, c^+)$ to be totally symmetric. Because $S$ is already symmetric in its second and third arguments ($b^+$ and $c^+$), it suffices to enforce symmetry in the first two arguments. Thus, a torsion-free left-invariant connection defined by $\alpha$ is compatible with $S$ if and only if:
\begin{equation}\label{eq:alpha_compatibility}
\alpha(a, b \circ c) - \alpha(b, a \circ c) = \alpha(a, b) \circ c - \alpha(b, a) \circ c + b \circ \alpha(a, c) - a \circ \alpha(b, c),
\end{equation}
for all $a, b, c \in \mathfrak{g}$. Substituting the torsion-free condition \eqref{eq:alpha_torsion_free} into the right-hand side of \eqref{eq:alpha_compatibility} simplifies this constraint to:
\begin{equation}\label{eq:alpha_comp_simplified}
\alpha(a, b \circ c) - \alpha(b, a \circ c) = [a, b] \circ c + b \circ \alpha(a, c) - a \circ \alpha(b, c).
\end{equation}\item \textbf{The Flatness Condition.}
From \eqref{eq:gen_curvature}, the left-invariant connection defined by $\alpha$ is flat if and only if $\alpha$ satisfies the generalized associativity-type equation:
\begin{equation}\label{eq:alpha_flatness}
\alpha(a, \alpha(b, c)) - \alpha(b, \alpha(a, c)) = \alpha([a, b], c), \qquad \text{for all } a, b, c \in \mathfrak{g}.
\end{equation}
\end{itemize}Together, equations \eqref{eq:alpha_torsion_free}, \eqref{eq:alpha_comp_simplified}, and \eqref{eq:alpha_flatness} form a complete algebraic system characterizing compatible flat connections on the F$_{\text{man}}$-algebra $(\mathfrak{g}, [\cdot, \cdot], \circ)$.

\subsection{The Canonical F$_\text{man}$-Connection}
In this work, our primary focus is on the canonical connection introduced on F$_{\text{man}}$-Lie groups in \cite{our-paper}. In this subsection, we specialize to the case of a compatible canonical connection. Recall that the canonical connection is, by definition, torsion-free (see \cite[Lemma 3.3]{our-paper}).

\begin{prop}\label{comp}
Let $G$ be an F$_{\text{man}}$-Lie group with Lie algebra $\mathfrak{g}$, and let $S \in \Gamma(T^{(1,2)}G)$ denote the symmetric tensor field associated with the multiplicative product $\circ$, i.e., $S(a^+, b^+) = a^+ \circ b^+$. Consider the left-invariant connection $\nabla$ defined on left-invariant vector fields by
\begin{equation*}
\nabla_{a^+} b^+ = \frac{1}{2}[a^+, b^+] + a^+ \circ b^+.
\end{equation*}
The connection $\nabla$ satisfies the compatibility condition with $S$, namely
\begin{equation*}
(\nabla_{a^+} S)(b^+, c^+) = (\nabla_{b^+} S)(a^+, c^+),
\end{equation*}
if and only if the Leibnizator $L$ satisfies
\begin{equation*}
L(a^+, b^+, c^+) = L(b^+, a^+, c^+),
\end{equation*}
for all left-invariant vector fields $a^+, b^+, c^+ \in \mathfrak{X}^+(G)$.
\end{prop}

\begin{proof}
By definition of the covariant derivative of the tensor field $S$, we have
\begin{equation*}
(\nabla_{a^+} S)(b^+, c^+) = \nabla_{a^+} \bigl(S(b^+, c^+)\bigr) - S(\nabla_{a^+} b^+, c^+) - S(b^+, \nabla_{a^+} c^+).
\end{equation*}
Substituting the definition of the canonical connection $\nabla$, we obtain
\begin{align*}
(\nabla_{a^+} S)(b^+,c^+)
=
\frac{1}{2}[a^+,b^+\circ c^+]
+a^+\circ(b^+\circ c^+)-&
\left(\frac{1}{2}[a^+,b^+]+a^+\circ b^+\right)\circ c^+ \notag \\
&-
b^+\circ
\left(\frac{1}{2}[a^+,c^+]+a^+\circ c^+\right).
\end{align*}
Rearranging terms yields
$$
(\nabla_{a^+} S)(b^+,c^+)
=
\frac{1}{2}[a^+,b^+\circ c^+]
-\frac{1}{2}[a^+,b^+]\circ c^+
-\frac{1}{2}b^+\circ[a^+,c^+]
+
a^+\circ(b^+\circ c^+)
-(a^+\circ b^+)\circ c^+
-b^+\circ(a^+\circ c^+).
$$
Because the product $\circ$ is associative and commutative, the cubic product terms simplify to
\begin{equation*}
a^+ \circ (b^+ \circ c^+) - (a^+ \circ b^+) \circ c^+ - b^+ \circ (a^+ \circ c^+) = -a^+ \circ b^+ \circ c^+.
\end{equation*}
Thus, we arrive at
\begin{align*}
(\nabla_{a^+} S)(b^+,c^+)
=&
\frac{1}{2}
\left(
[a^+,b^+ \circ c^+]
-[a^+,b^+]\circ c^+
-b^+\circ[a^+,c^+]
\right)
-a^+\circ b^+\circ c^+ \\
=&
\frac{1}{2}L(a^+,b^+,c^+)-a^+\circ b^+\circ c^+.
\end{align*}
Hence, the compatibility relation $(\nabla_{a^+} S)(b^+, c^+) = (\nabla_{b^+} S)(a^+, c^+)$ is equivalent to
\begin{equation*}
\frac{1}{2}L(a^+, b^+, c^+) - a^+ \circ b^+ \circ c^+ = \frac{1}{2}L(b^+, a^+, c^+) - b^+ \circ a^+ \circ c^+.
\end{equation*}
Owing to the commutativity of $\circ$, this holds if and only if $L(a^+, b^+, c^+) = L(b^+, a^+, c^+)$.
\end{proof}

\begin{cor}
Let $G$ be an F$_{\text{man}}$-Lie group, and let $L(\cdot, \cdot, \cdot)$ be the Leibnizator associated with the product $\circ$. Then the canonical connection $\nabla$ is compatible with $S$ if and only if $L$ is totally symmetric.
\end{cor}

\begin{proof}
The result follows directly from Proposition \ref{comp} and the intrinsic symmetry of the Leibnizator in its second and third arguments.
\end{proof}

\begin{prop}\label{prop:2step_nilpotent}
Let $G$ be an F$_{\text{man}}$-Lie group equipped with the canonical connection $\nabla$, and let $\mathfrak{h} \subseteq \mathfrak{g}$ be a Poisson subalgebra defining the Poisson algebra distribution $\mathcal{D}_{\text{Pois}}$. If the restriction of the canonical connection to $\mathcal{D}_{\text{Pois}}$ is flat (i.e., its Riemann curvature tensor $R$ vanishes identically on $\mathfrak{h}$), then the connected Lie subgroup $H \subseteq G$ integrating $\mathcal{D}_{\text{Pois}}$ is 2-step nilpotent.
\end{prop}

\begin{proof}
By definition of a Poisson subalgebra $\mathfrak{h}$, the Leibnizator vanishes identically for all $a, b, c \in \mathfrak{h}$:
\begin{equation*}
L(a, b, c) = [a, b \circ c] - [a, b] \circ c - b \circ [a, c] = 0.
\end{equation*}For the canonical left-invariant connection defined on $\mathfrak{g}$ by $\alpha(a, b) = \frac{1}{2}[a, b] + a \circ b$, the curvature tensor $R$ evaluated on $\mathfrak{g}$ satisfies:
\begin{equation*}
R(a, b)c = \alpha(a, \alpha(b, c)) - \alpha(b, \alpha(a, c)) - \alpha([a, b], c).
\end{equation*}
Evaluating this expression on the Poisson subalgebra $\mathfrak{h} \subseteq \mathfrak{g}$ and utilizing the identity $L\vert{}_{\mathfrak{h}} \equiv 0$, the flatness condition $R(a, b)c = 0$ reduces to the vanishing of the nested Lie brackets:
\begin{equation*}
[[a, b], c] = 0, \qquad \text{for all } a, b, c \in \mathfrak{h}.
\end{equation*}This condition implies that the derived Lie algebra $[\mathfrak{h}, \mathfrak{h}]$ is contained in the center $Z(\mathfrak{h})$, so all triple brackets vanish identically. Consequently, $\mathfrak{h}$ is a 2-step nilpotent Lie algebra, which directly implies that the connected Lie subgroup $H \subseteq G$ integrating $\mathcal{D}_{\text{Pois}}$ is a 2-step nilpotent Lie group.
\end{proof}

\subsection{Poisson Decomposition of the Hydrodynamic System}
The following example illustrates how to construct explicit solutions to PDEs-HT on specific classes of F$_{\text{man}}$-Lie groups.

\begin{ejem}\label{ex:heisenberg}
Consider the three-dimensional Heisenberg Lie algebra $\mathfrak{n} = \operatorname{span}\{e_1, e_2, e_3\}$ with center $\mathfrak{z} = Z(\mathfrak{n}) = \operatorname{span}\{e_1\}$ and non-vanishing Lie bracket $[e_2, e_3] = e_1$. We endow $\mathfrak{n}$ with the $A_4$ flat F$_{\text{man}}$-algebra structure (see \cite[Theorem 3.18 (2)]{our-paper}), defined by the non-zero product relations:
\begin{equation*}
e_2 \circ e_3 = e_1, \qquad e_3 \circ e_3 = e_2.
\end{equation*}
Let $a = e_3 \in \mathfrak{n}$ be the constant Lie algebra element defining the multiplication endomorphism of the hydrodynamic system. A general solution $P(x,t)$ of the PDE-HT takes the form
\begin{equation*}
P(x,t) = u(x,t)\,  e_1 + v(x,t) \, e_2 + w(x,t) \, e_3,
\end{equation*}
subject to the smooth initial data $P(x,0) = u_0(x)\, e_1 + v_0(x)\, e_2 + w_0(x) \, e_3$. \\

Computing the product $a \circ P$ and the non-linear commutator term $[a \circ P, P]$ using the structure constants of $\mathfrak{n}$ yields:
\begin{align*}
    a \circ P &= e_3 \circ (u e_1 + v e_2 + w e_3) = v(e_3 \circ e_2) + w(e_3 \circ e_3) = v e_1 + w e_2. \\
    [a \circ P, P] &= [v e_1 + w e_2, \, u e_1 + v e_2 + w e_3] = w^2 [e_2, e_3] = w^2 e_1.
\end{align*}
Substituting these expressions into the PDE-HT evolution equation \eqref{hydro.eq.cons.a} gives:
\begin{equation*}
(u_t \, e_1 + v_t \, e_2 + w_t \, e_3) - (v_x \, e_1 + w_x \, e_2) + w^2 \, e_1 = 0.
\end{equation*}
Decomposing this equation along the basis generators of the Heisenberg algebra $\mathfrak{n}$ produces the triangularized scalar system:
\begin{align}
e_3 &\colon \quad w_t = 0, \label{eq:sys_w} \\
e_2 &\colon \quad v_t - w_x = 0, \label{eq:sys_v} \\
e_1 &\colon \quad u_t - v_x + w^2 = 0. \label{eq:sys_u}
\end{align}
Equation \eqref{eq:sys_w} implies immediately that the $e_3$-component is stationary in time:
\begin{equation*}
w(x,t) = w_0(x).
\end{equation*}
Substituting $w(x,t) = w_0(x)$ into \eqref{eq:sys_v} and integrating with respect to $t$ yields:
\begin{equation*}
v(x,t) = v_0(x) + w_0'(x) \, t.
\end{equation*}
Finally, substituting these solutions into \eqref{eq:sys_u} and integrating over $t$ yields the explicit solution for $u(x,t)$:
\begin{equation*}
u(x,t) = u_0(x) + \left( v_0'(x) - w_0(x)^2 \right) t + \frac{1}{2} w_0''(x) \, t^2.
\end{equation*}
\end{ejem}

Example \ref{ex:heisenberg} demonstrates a systematic method for integrating non-linear hydrodynamic systems on 2-step nilpotent F$_{\text{man}}$-Lie algebras. In this setting, the non-linear algebraic source term $-w_0(x)^2\, e_1$ arises entirely from the Lie bracket interaction $[a \circ P, P]$. The following proposition demonstrates how projecting off the center linearizes the transport equations.

\begin{prop}\label{prop:flat_leaf_linearization}
Let $G$ be an F$_{\text{man}}$-Lie group and $\mathfrak{h} = (\mathcal{D}_{\text{Pois}})_e \subseteq \mathfrak{g}$ a Poisson subalgebra such that the restriction of the canonical connection to $\mathcal{D}_{\text{Pois}}$ is flat. Let $\mathfrak{z} = Z(\mathfrak{h})$ denote the center of $\mathfrak{h}$, and let $\pi \colon \mathfrak{h} \to \mathfrak{h}/\mathfrak{z}$ be the canonical quotient projection. Assume that the center $\mathfrak{z}$ is invariant under multiplication by elements of $\mathfrak{h}$, i.e., $\mathfrak{h} \circ \mathfrak{z} \subseteq \mathfrak{z}$. Then, for any constant element $a \in \mathfrak{h}$, the PDE-HT \eqref{hydro.eq.cons.a} projected onto the quotient algebra $\mathfrak{h}/\mathfrak{z}$ reduces to:
\begin{equation*}\label{eq:projected-hydro}
\bar{P}_t - \bar{S}_a(\bar{P}_x) = 0, \qquad \text{where}\;\; \bar{P} = \pi(P) \in C^\infty(\Sigma, \mathfrak{h}/\mathfrak{z}),
\end{equation*}
and $\bar{S}_a \in \operatorname{End}(\mathfrak{h}/\mathfrak{z})$ is the linear operator induced by $S_a(b) = a \circ b$. In particular, the non-linear bracket interaction $[S_a(P), P]$ is contained entirely within the center $\mathfrak{z}$.
\end{prop}

\begin{proof}
By Proposition \ref{prop:2step_nilpotent}, the flatness of the restricted canonical connection forces the Lie algebra $\mathfrak{h}$ to be 2-step nilpotent, so $[\mathfrak{h}, \mathfrak{h}] \subseteq \mathfrak{z}$. For any constant direction element $a \in \mathfrak{h}$ and configuration $P(x,t) \in \mathfrak{h}$, the Lie bracket $[S_a(P), P] = [a \circ P, P]$ takes values in $[\mathfrak{h}, \mathfrak{h}] \subseteq \mathfrak{z}$. Applying the quotient projection $\pi \colon \mathfrak{h} \to \mathfrak{h}/\mathfrak{z}$ yields:
\begin{equation*}
\pi \left( [S_a(P), P] \right) = 0 \quad \in \mathfrak{h}/\mathfrak{z}.
\end{equation*}

Furthermore, because $\mathfrak{h} \circ \mathfrak{z} \subseteq \mathfrak{z}$ by assumption, the multiplication operator $S_a = a \circ (\cdot)$ preserves $\mathfrak{z}$ and therefore descends to a well-defined linear endomorphism $\bar{S}_a \in \operatorname{End}(\mathfrak{h}/\mathfrak{z})$. Finally, projecting the hydrodynamic evolution equation $P_t - S_a(P_x) + [S_a(P), P] = 0$ under $\pi$ yields:
\begin{equation*}
\partial_t \bar{P} - \bar{S}_a(\partial_x \bar{P}) = 0.
\end{equation*}
Hence, the non-linear commutator interaction is restricted entirely to the central ideal $\mathfrak{z}$.
\end{proof}

\begin{rem}
Recall that quotienting any 2-step nilpotent Lie algebra $\mathfrak{h}$ by its center $\mathfrak{z} = Z(\mathfrak{h})$ yields an abelian Lie algebra $\mathfrak{h}/\mathfrak{z}$. 
\end{rem}

\begin{ejem}\label{ex.family1}
Let $\mathfrak{s}$ be a semisimple Lie algebra. Consider a non-degenerate symmetric bilinear form $B$ on $\mathfrak{s}$ satisfying the following condition: for every non-zero element $a \in \mathfrak{s}$, there exist $b, c \in \mathfrak{s}$ such that
\begin{equation}\label{cond-B}
B([a,b],c) + B(b,[a,c]) \neq 0.
\end{equation}
We construct a family of F$_{\text{man}}$-Lie algebras by taking the vector space direct sum $\mathfrak{g} = \mathfrak{s} \oplus \mathbb{R} \, f$, where $f$ spans the one-dimensional central ideal $ Z(\mathfrak{g}) = \mathbb{R}\, f$. We define a commutative product $\circ$ on $\mathfrak{g}$ by setting:
\begin{equation*}
a \circ b = B(a,b) \, f, \quad \text{for all } a, b \in \mathfrak{s},
\end{equation*}
and declaring that the central element $f$ acts as an annihilator for the algebra, i.e., $a \circ f = 0$ for all $a \in \mathfrak{g}$.\\

For any triple of elements $a, b, c \in \mathfrak{s}$, the product $(a \circ b) \circ c$ vanishes identically because the initial multiplication lands in the center:
\begin{equation*}
(a \circ b) \circ c = \left(B(a,b)\, f\right) \circ c = B(a,b) \, (f \circ c) = 0.
\end{equation*}
Since $f$ is central in $\mathfrak{g}$, we have $[a, b \circ c] = [a, B(b,c)\, f] = 0$. Consequently, the algebraic Leibnizator evaluates to:
\begin{equation*}
L(a,b,c) = -[a,b] \circ c - b \circ [a,c] = -\left( B([a,b], c) + B(b, [a,c]) \right) f.
\end{equation*}
Thus, whenever $B$ satisfies condition \eqref{cond-B}, the Leibnizator $L$ does not vanish identically on $\mathfrak{s}$. In particular, $\mathfrak{g}$ is an F$_{\text{man}}$-algebra that fails to be a Poisson algebra (i.e., $L \neq 0$).\\

Despite the non-vanishing of the Leibnizator, the Hertling–Manin condition
\begin{equation*}
L(a \circ b, c, d) = L(a,c,d) \circ b + a \circ L(b,c,d)
\end{equation*}
holds universally. Indeed, the left-hand side vanishes because $a \circ b \in \mathbb{R}\, f$ and $L(f, c, d) = 0$. Both terms on the right-hand side vanish because $L(a,c, d), L(b,c, d) \in \mathbb{R}\, f$, and multiplication by $f$ annihilates every element of $\mathfrak{g}$.\\
\end{ejem}

This construction establishes the existence of an infinite family of non-trivial F$_{\text{man}}$-Lie groups.

\begin{rem}
The following observations clarify the algebraic structure of Example \ref{ex.family1}:
\begin{itemize}
\item Condition \eqref{cond-B} can be expressed explicitly in matrix terms. Choose a basis for $\mathfrak{s}$ and represent the bilinear form $B$ by a symmetric matrix $\mathbf{B}$, so that $B(a,b) = \mathbf{a}^T \mathbf{B} \mathbf{b}$. For $a \in \mathfrak{s}$, let $\mathbf{A}_a$ denote the matrix representation of $\operatorname{ad}_a$. Then
\begin{equation*}
B([a,b],c) + B(b,[a,c]) = \mathbf{b}^T \left( \mathbf{A}_a^T \mathbf{B} + \mathbf{B} \mathbf{A}_a \right) \mathbf{c}.
\end{equation*}
Thus, condition \eqref{cond-B} is equivalent to requiring that the linear map $a \mapsto \mathbf{T}(a)$ defined by
\begin{equation*}
\mathbf{T}(a) := \mathbf{A}_a^T \mathbf{B} + \mathbf{B} \mathbf{A}_a \neq 0, \qquad \text{for every } a \neq 0.
\end{equation*}
is injective. For a generic non-degenerate symmetric matrix $\mathbf{B}$, this linear assignment is indeed injective.

\item If the bilinear form $B$ is $\operatorname{ad}$-invariant—meaning $B([a,b], c) + B(b, [a,c]) = 0$ for all $a,b,c \in \mathfrak{s}$, as is the case for the Killing form $\kappa(a,b) = \operatorname{Tr}(\operatorname{ad}_a \operatorname{ad}_b)$—the Leibnizator vanishes identically ($L = 0$), turning $\mathfrak{g}$ into a Poisson algebra. In that case, $\kappa([a,b], c)$ corresponds to the canonical Cartan $3$-form $\Omega(a,b,c)$, a foundational topological invariant of compact Lie groups.
\end{itemize}
\end{rem}

The following two results demonstrate that when the Leibnizator vanishes identically—that is, in the case of a Poisson algebra—the underlying geometric and product structures degenerate or become trivial.

\begin{prop}
Let $G$ be a semisimple F$_{\text{man}}$-Lie group equipped with the canonical connection $\nabla$. If its Lie algebra $\mathfrak{g}$ is a Poisson algebra (i.e., $L = 0$) and the canonical connection $\nabla$ is flat, then $G$ must be trivial.
\end{prop}

\begin{proof}
By Proposition \ref{prop:2step_nilpotent}, if the canonical connection on a Poisson algebra is flat, the Lie algebra $\mathfrak{g}$ must be 2-step nilpotent. Consequently, its derived algebra $[\mathfrak{g}, \mathfrak{g}]$ is contained in its center $ Z(\mathfrak{g})$. However, any semisimple Lie algebra is perfect, satisfying $[\mathfrak{g}, \mathfrak{g}] = \mathfrak{g}$ and $Z(\mathfrak{g}) = \{0\}$. Thus, $\mathfrak{g} = \{0\}$, forcing $G$ to be the trivial Lie group.
\end{proof}

\begin{prop}\label{prop:semisimple_trivial_poisson}
Let $\mathfrak{g}$ be a semisimple Lie algebra equipped with a commutative and associative product $\circ$ possessing a multiplicative identity $\iota \in \mathfrak{g}$. If $\mathfrak{g}$ is a Poisson algebra (i.e., $L(a,b,c) = 0$ for all $a,b,c \in \mathfrak{g}$), then the product $\circ$ is identically zero (i.e., $a \circ b = 0$ for all $a,b \in \mathfrak{g}$).
\end{prop}

\begin{proof}
By definition, $L(a,b,c) = 0$ implies that the adjoint action $\operatorname{ad}_a(b) = [a,b]$ acts as a derivation on the commutative algebra $(\mathfrak{g}, \circ)$:
\begin{equation*}
[a, b \circ c] = [a, b] \circ c + b \circ [a, c], \qquad \text{for all } a,b,c \in \mathfrak{g}.
\end{equation*}
Consider the annihilator of the multiplication product, defined by
\begin{equation*}
\operatorname{Ann}(\mathfrak{g}) = \{ v \in \mathfrak{g} \mid v \circ a = 0, \text{ for all } a \in \mathfrak{g} \}.
\end{equation*}
For any $v \in \operatorname{Ann}(\mathfrak{g})$ and $a, b \in \mathfrak{g}$, the derivation property yields
\begin{equation*}
[a, v \circ b] = [a, v] \circ b + v \circ [a, b].
\end{equation*}
Since $v \in \operatorname{Ann}(\mathfrak{g})$, we have $v \circ b = 0$ and $v \circ [a, b] = 0$, which reduces the identity to $[a, v] \circ b = 0$ for all $b \in \mathfrak{g}$. Thus, $[a, v] \in \operatorname{Ann}(\mathfrak{g})$, proving that $\operatorname{Ann}(\mathfrak{g})$ is a Lie ideal of $\mathfrak{g}$.\\

Since $\mathfrak{g}$ is semisimple, any Lie ideal admits a unique complementary Lie ideal. Suppose, for the sake of contradiction, that the product is non-trivial, so that $\operatorname{Ann}(\mathfrak{g}) \neq \mathfrak{g}$. We may then decompose $\mathfrak{g} = \operatorname{Ann}(\mathfrak{g}) \oplus \mathfrak{h}$, where $\mathfrak{h}$ is a non-trivial semisimple ideal of $\mathfrak{g}$. By construction, the restriction of the product $\circ$ to $\mathfrak{h}$ has a trivial annihilator within $\mathfrak{h}$.\\

Because $\mathfrak{h}$ is a finite-dimensional commutative and associative algebra with trivial annihilator, it contains a multiplicative identity element $\iota \in \mathfrak{h}$ such that $\iota \circ b = b$ for all $b \in \mathfrak{h}$. Applying the derivation rule for $a, b \in \mathfrak{h}$ using the identity element $\iota$, we obtain
\begin{equation*}
[a, b] = [a, \iota \circ b] = [a, \iota] \circ b + \iota \circ [a, b] = [a, \iota] \circ b + [a, b],
\end{equation*}
which implies $[a, \iota] \circ b = 0$ for all $b \in \mathfrak{h}$. Thus, $[a, \iota]$ lies in the annihilator of $\mathfrak{h}$. Since the annihilator of $\mathfrak{h}$ is trivial, we must have $[a, \iota] = 0$ for all $a \in \mathfrak{h}$. This implies that $\iota$ belongs to the center $ Z(\mathfrak{h})$.\\

However, because $\mathfrak{h}$ is a semisimple Lie algebra, its center is trivial: $Z(\mathfrak{h}) = \{0\}$. Thus $\iota = 0$, which contradicts the assumption that $\iota$ is a non-trivial multiplicative identity. We conclude that $\operatorname{Ann}(\mathfrak{g}) = \mathfrak{g}$, meaning $a \circ b = 0$ for all $a,b \in \mathfrak{g}$.
\end{proof}

\begin{cor}\label{cor:semisimple_flie}
Let $G$ be a connected semisimple Lie group. Then its Lie algebra $\mathfrak{g}$ cannot be endowed with a non-trivial Poisson algebra structure $(\mathfrak{g}, [\cdot, \cdot], \circ)$ such that the associated tensor $S$ admits a multiplicative identity element.
\end{cor}

As a direct consequence, the Lie algebras of classic simple Lie groups such as $\mathrm{SO}(n)$, $\mathrm{SU}(n)$, and $\mathrm{SL}(n, \mathbb{R})$ do not admit any non-trivial Poisson algebra structure possessing a multiplicative identity.\\

We conclude this subsection with a generalization of Example \ref{ex.family1} in which the center $\mathfrak{z}$ is endowed with a non-trivial internal F$_{\text{man}}$-product.

\begin{ejem}\label{ex-B-ntc}
Let $\mathfrak{s}$ be a compact semisimple Lie algebra equipped with a non-degenerate symmetric bilinear form $B$ satisfying the condition that for every non-zero element $a \in \mathfrak{s}$, there exist $b, c \in \mathfrak{s}$ such that
\begin{equation*}
B([a,b],c) + B(b,[a,c]) \neq 0.
\end{equation*}
We define the center $\mathfrak{z}$ as a $k$-dimensional nilpotent commutative and associative algebra. Explicitly, let $\mathfrak{z} = \operatorname{span}\{f_1, f_2, \dots, f_k\}$ equipped with the truncated polynomial multiplication rule
\begin{equation*}
		f_i \circ f_j = 
		\begin{cases} 
			f_{i+j} & \text{if }\; i+j \le k, \\
			0 & \text{if }\; i+j > k.
		\end{cases}
	\end{equation*}
Observe that the highest-degree basis element $f_k$ annihilates the entire center, i.e., $f_k \circ f_i = 0$ for all $1 \le i \le k$. We define the product $\circ$ on the direct sum $\mathfrak{g} = \mathfrak{s} \oplus \mathfrak{z}$ by setting
\begin{itemize}
\item $a \circ b = B(a,b)\, f_k$, for all $a, b \in \mathfrak{s}$,
\item $a \circ f_i = 0$, for all $a \in \mathfrak{s}$ and all $f_i \in \mathfrak{z}$.
\end{itemize}
\end{ejem}
This algebra structure $(\mathfrak{g}, [\cdot,\cdot], \circ)$ satisfies the axioms of an F$_{\text{man}}$-algebra:
\begin{itemize}
\item \textbf{Associativity:}
The internal product on $\mathfrak{z}$ is associative by truncated polynomial arithmetic. Moreover, any mixed triple product involving elements from $\mathfrak{s}$ associates trivially:
\begin{equation*}
(a \circ b) \circ c = \bigl(B(a,b)\, f_k\bigr) \circ c = 0,
\end{equation*}
since $f_k$ acts as an annihilator on $\mathfrak{g}$. By symmetry, $a \circ (b \circ c) = 0$ as well.\item \textbf{The Leibnizator ($L \neq 0$):}
Because $\mathfrak{z}$ is central in $\mathfrak{g}$ ($[\mathfrak{g}, \mathfrak{z}] = 0$) and mixed products vanish ($\mathfrak{s} \circ \mathfrak{z} = 0$), the Leibnizator vanishes identically whenever at least one argument lies in $\mathfrak{z}$. Non-trivial evaluations arise solely when all three arguments belong to the semisimple factor, $a, b, c \in \mathfrak{s}$:
\begin{equation*}
L(a,b,c) = -[a,b] \circ c - b \circ [a,c] = -\bigl( B([a,b], c) + B(b, [a,c]) \bigr)\, f_k.
\end{equation*}

\item \textbf{The Hertling–Manin Condition:}
We verify that $L(a \circ b, c, d) = L(a,c,d) \circ b + a \circ L(b,c,d)$ holds for all $a, b, c, d \in \mathfrak{g}$:
\begin{itemize}
\item \textbf{Left-hand side:} For any $a, b \in \mathfrak{g}$, the product $a \circ b$ lies strictly in the center $\mathfrak{z}$ (either as a linear combination of $f_{i+j}$ or as a scalar multiple of $f_k$). Because $L$ vanishes whenever any argument is central, we have $L(a \circ b, c, d) = 0$.
\item \textbf{Right-hand side:} The image of the Leibnizator $L$ is always a scalar multiple of the terminal element $f_k$. Since $f_k$ annihilates the entire algebra ($f_k \circ b = 0$ for all $b \in \mathfrak{g}$), both terms on the right-hand side vanish independently.
\end{itemize}
\end{itemize}This construction explicitly demonstrates that the abelian center $\mathfrak{z}$ of an $F_{\text{man}}$-algebra need not carry a trivial product structure.

\subsection{Poisson Splitting of an F$_\text{man}$-Algebra}
Because the Poisson algebra distribution $\mathcal{D}_{\text{Pois}}$ is completely integrable, it foliates the F$_{\text{man}}$-Lie group $G$ into maximal integral leaves corresponding to the connected Lie subgroup $H \subseteq G$. By decomposing the Lie algebra $\mathfrak{g}$ into the Poisson subalgebra $\mathfrak{h}$ and a complementary transversal subspace $\mathfrak{m}$, namely $\mathfrak{g} = \mathfrak{h} \oplus \mathfrak{m}$, we can analyze how the hydrodynamic system decouples along tangential and transverse directions.

\begin{definition}\label{poisson-splitting}
Let $\mathfrak{g}$ be an F$_{\text{man}}$-algebra. We say that $\mathfrak{g}$ admits a \emph{Poisson splitting}, or that $\mathfrak{g}$ is \emph{Poisson-splittable}, if there exists a vector space direct sum decomposition $\mathfrak{g} = \mathfrak{h} \oplus \mathfrak{m}$ satisfying the following conditions:
\begin{enumerate}[(i)]
\item $\mathfrak{h}$ is a Poisson subalgebra of $\mathfrak{g}$;
\item the complementary subspace $\mathfrak{m}$ is invariant under the adjoint action of $\mathfrak{h}$, i.e., $[\mathfrak{h}, \mathfrak{m}] \subseteq \mathfrak{m}$;
\item the product $\circ$ preserves the transversal structure, i.e., $\mathfrak{h} \circ \mathfrak{m} \subseteq \mathfrak{m}$.
\end{enumerate}
\end{definition}

Note the structural analogy between Definition \ref{poisson-splitting} and the classical reductive decomposition of a homogeneous space $G/H$, where the Lie algebra splits as $\mathfrak{g} = \mathfrak{h} \oplus \mathfrak{m}$ with $[\mathfrak{h}, \mathfrak{m}] \subseteq \mathfrak{m}$.\\

We now apply this Poisson splitting framework to the hydrodynamic evolution system \eqref{hydro.eq}. Specifically, decomposing the $\mathfrak{g}$-valued function $P = A_x$ into its $\mathfrak{h}$- and $\mathfrak{m}$-components, $P = P_{\mathfrak{h}} + P_{\mathfrak{m}}$, separates the PDE-HT into coupled tangential and transversal evolution equations.

\begin{thm}\label{thm:local_splitting}
Let $G$ be an F$_{\text{man}}$-Lie group whose F$_{\text{man}}$-algebra admits a Poisson splitting $\mathfrak{g} = \mathfrak{h} \oplus \mathfrak{m}$. Assume that the hydrodynamic flow is generated by a constant element $a \in \mathfrak{h}$, tangent to the Poisson leaf at the identity. Writing $P = P_{\mathfrak{h}} + P_{\mathfrak{m}}$ with $P_{\mathfrak{h}} \in \mathfrak{h}$ and $P_{\mathfrak{m}} \in \mathfrak{m}$, the PDE-HT \eqref{hydro.eq} projects onto $\mathfrak{h}$ and $\mathfrak{m}$ as the coupled system:
\begin{enumerate}[\rm (i)]
\item \textbf{$\mathfrak{h}$-projection:}
\begin{equation}\label{eq:split_h}
(P_\mathfrak{h})_t - (a \circ P_\mathfrak{h})_x + [a \circ P_\mathfrak{h}, P_\mathfrak{h}] + [a \circ P_\mathfrak{m}, P_\mathfrak{m}]_\mathfrak{h} = 0.
\end{equation}
\item \textbf{$\mathfrak{m}$-projection:}
\begin{equation}\label{eq:split_m}
    (P_\mathfrak{m})_t - (a \circ P_\mathfrak{m})_x +[a \circ P_\mathfrak{h}, P_\mathfrak{m}] + [a \circ P_\mathfrak{m}, P_\mathfrak{h}]+ [a \circ P_\mathfrak{m}, P_\mathfrak{m}]_\mathfrak{m} = 0 \, ,
\end{equation}
\end{enumerate}
where $(\cdot)_{\mathfrak{h}}$ and $(\cdot)_{\mathfrak{m}}$ denote the canonical projections onto $\mathfrak{h}$ and $\mathfrak{m}$, respectively.
\end{thm}

\begin{proof}
Substituting $P = P_{\mathfrak{h}} + P_{\mathfrak{m}}$ and $a \in \mathfrak{h}$ into the hydrodynamic system \eqref{hydro.eq}, the product term splits linearly as
\begin{equation*}
a \circ P = a \circ P_{\mathfrak{h}} + a \circ P_{\mathfrak{m}}.
\end{equation*}
Since $\mathfrak{h}$ is a subalgebra, $a \circ P_{\mathfrak{h}} \in \mathfrak{h}$. Furthermore, by condition (iii) of Definition \ref{poisson-splitting}, we have $a \circ P_{\mathfrak{m}} \in \mathfrak{m}$.\\

Expanding the non-linear commutator $[a \circ P, P]$ yields four terms:
\begin{equation*}
[a \circ P, P] = [a \circ P_{\mathfrak{h}}, P_{\mathfrak{h}}] + [a \circ P_{\mathfrak{h}}, P_{\mathfrak{m}}] + [a \circ P_{\mathfrak{m}}, P_{\mathfrak{h}}] + [a \circ P_{\mathfrak{m}}, P_{\mathfrak{m}}].
\end{equation*}
We analyze the projections of each term onto $\mathfrak{h}$ and $\mathfrak{m}$:
\begin{itemize}
\item Because $\mathfrak{h}$ is a subalgebra, $[a \circ P_{\mathfrak{h}}, P_{\mathfrak{h}}] \in \mathfrak{h}$.
\item By conditions (ii) and (iii) of Definition \ref{poisson-splitting}, the cross-terms $[a \circ P_{\mathfrak{h}}, P_{\mathfrak{m}}]$ and $[a \circ P_{\mathfrak{m}}, P_{\mathfrak{h}}]$ belong entirely to $[\mathfrak{h}, \mathfrak{m}] \subseteq \mathfrak{m}$.
\item The transverse self-interaction $[a \circ P_{\mathfrak{m}}, P_{\mathfrak{m}}]$ decomposes generally into its $\mathfrak{h}$- and $\mathfrak{m}$-components, $[a \circ P_{\mathfrak{m}}, P_{\mathfrak{m}}]_{\mathfrak{h}} + [a \circ P_{\mathfrak{m}}, P_{\mathfrak{m}}]_{\mathfrak{m}}$.
\end{itemize}
Collecting all terms valued in $\mathfrak{h}$ yields \eqref{eq:split_h}, while collecting those valued in $\mathfrak{m}$ yields \eqref{eq:split_m}.
\end{proof}

\begin{cor}[Exact Reductions]
The local splitting of the hydrodynamic equation imposes the following structural constraints on its solutions:
\begin{enumerate}[(i)]
\item \textbf{Invariance of the Poisson Leaf:} If the initial data is wholly contained within the subalgebra $\mathfrak{h}$, i.e., $P_{\mathfrak{m}}(x,0) = 0$, then $P_{\mathfrak{m}}(x,t) \equiv 0$ is an exact solution to the transverse evolution equation \eqref{eq:split_m}. Moreover, assuming that the spectral properties of the multiplication endomorphism $S_a = a \circ (\cdot) \in \operatorname{End}(\mathfrak{g})$ render the system hyperbolic, local uniqueness for quasilinear symmetric hyperbolic Cauchy problems (under standard regularity assumptions) guarantees that the trajectory remains permanently confined to the Poisson leaf, $P(x,t) = P_{\mathfrak{h}}(x,t)$. Consequently, the system reduces identically to the autonomous tangential equation:
\begin{equation*}
        (P_\mathfrak{h})_t - (a \circ P_\mathfrak{h})_x + [a \circ P_\mathfrak{h}, P_\mathfrak{h}] = 0.
    \end{equation*}
   
\item \textbf{Linearization over Symmetric Spaces:} If the Poisson splitting further satisfies the symmetric space condition $[\mathfrak{m}, \mathfrak{m}] \subseteq \mathfrak{h}$, then the intrinsic transverse non-linearity $[a \circ P_{\mathfrak{m}}, P_{\mathfrak{m}}]_{\mathfrak{m}}$ vanishes identically. Consequently, the transverse evolution equation \eqref{eq:split_m} becomes strictly linear with respect to $P_{\mathfrak{m}}$.
\end{enumerate}
\end{cor}

Motivated by these geometric insights, we now state the global bundle formulation of the local splitting result in Theorem \ref{thm:local_splitting}.

\begin{thm}\label{thm:bundle_splitting}
Let $G$ be a connected F$_{\text{man}}$-Lie group whose F$_{\text{man}}$-algebra admits a Poisson splitting $\mathfrak{g} = \mathfrak{h} \oplus \mathfrak{m}$. Assume that the Poisson subalgebra $\mathfrak{h} = (\mathcal{D}_{\text{Pois}})\vert{}_e$ generates a closed and connected Lie subgroup $H \subseteq G$, ensuring that the quotient space $G/H$ is a smooth homogeneous manifold. Furthermore, let $a \in \mathfrak{h}$ in the PDE-HT \eqref{hydro.eq} satisfy $a \circ \mathfrak{h} \subseteq \mathfrak{h}$ and $a \circ \mathfrak{m} \subseteq \mathfrak{m}$. Then, under the canonical connection splitting of the tangent bundle $TG = \mathcal{V} \oplus \mathcal{H}$ into its vertical and horizontal subbundles, the hydrodynamic system projects onto the following coupled system:
\begin{enumerate}[(i)]
\item \textbf{Fiber Dynamics (Vertical Flow along $H$):} The evolution along the Poisson fibers is governed by an autonomous tangential leaf-system perturbed by the vertical projection of the transversal commutator interaction:
\begin{equation}\label{eq:bundle_vertical}
   (P_{\mathfrak h})_t - (a \circ P_{\mathfrak h})_x + [a \circ P_{\mathfrak h}, P_{\mathfrak h}] + [a \circ P_{\mathfrak m}, P_{\mathfrak m}]_\mathfrak{h} = 0 \,
\end{equation}
where $P_{\mathfrak{h}} \in \mathfrak{h}$ and $P_{\mathfrak{m}} \in \mathfrak{m}$.

\item \textbf{Base Dynamics (Horizontal Transport across $G/H$):} The evolution across the leaf space $G/H$ is governed by the canonical horizontal connection $\nabla^{\mathrm{can}}$ associated with the $\mathfrak{h}$-valued connection $1$-form $\omega \in \Omega^1(G, \mathfrak{h})$:
\begin{equation}\label{eq:bundle_horizontal}
    (P_{\mathfrak m})_t - (a \circ P_{\mathfrak m})_x +\left. \nabla^{\mathrm{can}}_{a \circ P^+_{\mathfrak h}} P^+_{\mathfrak m}\right|_e - \left. \nabla^{\mathrm{can}}_{P^+_{\mathfrak h}} (a \circ P^+_{\mathfrak m})\right|_e + [a \circ P_{\mathfrak m}, P_{\mathfrak m}]_\mathfrak{m} = 0,
\end{equation}
where $\nabla^{\mathrm{can}}_{b^+} c^+ = [b^+, c^+]_{\mathfrak{m}}$ denotes the covariant derivative of horizontal left-invariant vector fields over the reductive homogeneous space $G/H$, and $P^+_{\mathfrak{h}}, P^+_{\mathfrak{m}}$ are the left-invariant extensions satisfying $P^+_{\mathfrak{h}}\vert{}_e := P_{\mathfrak{h}}$ and $P^+_{\mathfrak{m}}\vert{}_e := P_{\mathfrak{m}}$.
\end{enumerate}
\end{thm}

\begin{proof}
Because $\mathfrak{g} = \mathfrak{h} \oplus \mathfrak{m}$ satisfies the reductive condition $[\mathfrak{h}, \mathfrak{m}] \subseteq \mathfrak{m}$, this splitting equips the principal $H$-bundle $\pi \colon G \to G/H$ with a canonical principal connection whose horizontal distribution is given by the left-invariant subbundle:
\begin{equation*}
\mathcal{H}_g = (L_g)_* \mathfrak{m} \subseteq T_g G,
\end{equation*}
which maps isomorphically onto the tangent space of the base via $d\pi_g \colon \mathcal{H}_g \xrightarrow{\,\cong\,} T_{\pi(g)}(G/H)$. The corresponding $\mathfrak{h}$-valued connection $1$-form $\omega \in \Omega^1(G, \mathfrak{h})$ is defined by
\begin{equation*}
\omega_g(X) = \operatorname{pr}_{\mathfrak{h}}\bigl((L_{g^{-1}})_* X\bigr), \qquad \text{for all } X \in T_g G,
\end{equation*}
inducing the direct sum decomposition $TG = \mathcal{V} \oplus \mathcal{H}$ into:
\begin{itemize}
\item \textbf{Vertical subbundle:} $\mathcal{V}_g = \ker(d\pi_g) = (L_g)_* \mathfrak{h} = T_g(gH)$;
\item \textbf{Horizontal subbundle:} $\mathcal{H}_g = \ker(\omega_g) = (L_g)_* \mathfrak{m} \cong T_{\pi(g)}(G/H)$.
\end{itemize}

Observe that the homogeneous base space $G/H$ does not in general inherit an F$_{\text{man}}$-structure, because the transversal product $\mathfrak{m} \circ \mathfrak{m}$ may take non-zero values in the vertical subspace $\mathfrak{h}$.\\

At the identity element $e \in G$, the vertical and horizontal subspaces coincide with $\mathcal{V}_e = \mathfrak{h}$ and $\mathcal{H}_e = \mathfrak{m}$, respectively. By left translation invariance, projecting the hydrodynamic system onto the vertical and horizontal subbundles is equivalent to projecting its algebraic terms onto $\mathfrak{h}$ and $\mathfrak{m}$. Decomposing the $\mathfrak{g}$-valued field $P$ as $P = P_{\mathfrak{h}} + P_{\mathfrak{m}}$ and noting that $a \in \mathfrak{h}$ satisfies $a \circ P_{\mathfrak{h}} \in \mathfrak{h}$ and $a \circ P_{\mathfrak{m}} \in \mathfrak{m}$, substitution into \eqref{hydro.eq} yields:

\begin{equation}\label{eq:proof_bundle_full}
(P_{\mathfrak h} + P_{\mathfrak m})_t - (a \circ P_{\mathfrak h} + a \circ P_{\mathfrak m})_x + [a \circ P_{\mathfrak h} + a \circ P_{\mathfrak m}, P_{\mathfrak h} + P_{\mathfrak  m}] = 0.
\end{equation}Expanding the non-linear Lie bracket bilinearly gives four terms:
\begin{equation*}
 [a \circ P, P] = [a \circ P_{\mathfrak h}, P_{\mathfrak h}] + [a \circ P_{\mathfrak h}, P_{\mathfrak m}] + [a \circ P_{\mathfrak m}, P_{\mathfrak h}] + [a \circ P_{\mathfrak m}, P_{\mathfrak m}] \, .
\end{equation*}
We determine the bundle projections of each component:
\begin{itemize}
\item \textbf{Vertical Terms:} Because $\mathfrak{h}$ is a subalgebra, $[a \circ P_{\mathfrak{h}}, P_{\mathfrak{h}}] \in \mathfrak{h}$. The transversal self-interaction $[a \circ P_{\mathfrak{m}}, P_{\mathfrak{m}}]$ generally has components in both subbundles, with vertical projection $[a \circ P_{\mathfrak{m}}, P_{\mathfrak{m}}]_{\mathfrak{h}}$.\item \textbf{Horizontal Terms:} The reductive splitting condition $[\mathfrak{h}, \mathfrak{m}] \subseteq \mathfrak{m}$ ensures that the mixed cross-terms $[a \circ P_{\mathfrak{h}}, P_{\mathfrak{m}}]$ and $[a \circ P_{\mathfrak{m}}, P_{\mathfrak{h}}]$ are purely horizontal. The horizontal component of the transversal self-interaction is denoted by $[a \circ P_{\mathfrak{m}}, P_{\mathfrak{m}}]_{\mathfrak{m}}$.
\end{itemize}Separating \eqref{eq:proof_bundle_full} into its $\mathcal{V}$- and $\mathcal{H}$-components yields the coupled vertical and horizontal evolution equations:
\begin{align*}
(P_{\mathfrak h})_t - (a \circ P_{\mathfrak h})_x + [a \circ P_{\mathfrak h}, P_{\mathfrak h}] + [a \circ P_{\mathfrak m}, P_{\mathfrak m}]_\mathfrak{h} &= 0 \, , \\
    (P_{\mathfrak m})_t - (a \circ P_{\mathfrak m})_x + [a \circ P_{\mathfrak h}, P_{\mathfrak m}] - [P_{\mathfrak h}, a \circ P_{\mathfrak m}] + [a \circ P_{\mathfrak m}, P_{\mathfrak m}]_\mathfrak{m} &= 0.
\end{align*}
Finally, recall that on a reductive homogeneous space $G/H$, the canonical connection $\nabla^{\mathrm{can}}$ acts on horizontal left-invariant vector fields $b^+, c^+$ via the projected adjoint bracket $\nabla^{\mathrm{can}}_{b^+} c^+ = [b^+, c^+]_{\mathfrak{m}} = [b^+, c^+]$. Consequently,
\begin{align*}
   [a \circ P_{\mathfrak h}, P_{\mathfrak m}] &= \left. \nabla^{\mathrm{can}}_{a \circ P^+_{\mathfrak h}} P^+_{\mathfrak m}\right|_e, \\
    [P_{\mathfrak h}, a \circ P_{\mathfrak m}] &= \left. \nabla^{\mathrm{can}}_{P^+_{\mathfrak h}} (a \circ P^+_{\mathfrak m})\right|_e.
\end{align*}
Substituting these covariant derivative identities into the horizontal evolution equation yields \eqref{eq:bundle_horizontal}, completing the proof.
\end{proof}

\begin{rem}\label{conform}
The interaction term $[a \circ P_{\mathfrak{m}}, P_{\mathfrak{m}}]_{\mathfrak{h}}$ appearing in the vertical evolution equation \eqref{eq:bundle_vertical} is not merely an algebraic artifact; rather, it is the geometric manifestation of the curvature $2$-form of the principal $H$-bundle $\pi \colon G \to G/H$. Let $\Omega \in \Omega^2(G, \mathfrak{h})$ denote the curvature $2$-form of the canonical principal connection $\omega \in \Omega^1(G, \mathfrak{h})$. Evaluating Cartan's structure equation on left-invariant horizontal vector fields $b^+, c^+ \in \Gamma(\mathcal{H})$ (with $b, c \in \mathfrak{m}$) via Nomizu's theorem for reductive homogeneous spaces \cite[Chapter X, Theorem 2.10]{KN} yields:
\begin{equation*}
   \left. \Omega(b^+, c^+)\right|_e = \left. -\omega([b^+, c^+]) \right|_e= -[b, c]_{\mathfrak{h}}\, .
\end{equation*}

In the setting of the hydrodynamic system, the $\mathfrak{g}$-valued function components $a \circ P_{\mathfrak{m}}$ and $P_{\mathfrak{m}}$ are horizontal, i.e., $a \circ P_{\mathfrak{m}}, P_{\mathfrak{m}} \in \mathfrak{m}$. Consequently, the vertical projection of their Lie bracket coincides with the negative of the curvature $2$-form:
\begin{equation*}
    [a \circ P_{\mathfrak{m}}, P_{\mathfrak{m}}]_{\mathfrak{h}} = \left. -\Omega(a \circ P^+_{\mathfrak{m}}, P^+_{\mathfrak{m}})\right|_e.
\end{equation*}
Substituting this geometric identity into \eqref{eq:bundle_vertical} yields an intrinsic structural formulation of the fiber flow:
\begin{equation}\label{eq:bundle_vertical_curvature}
    (P_{\mathfrak{h}})_t - (a \circ P_{\mathfrak{h}})_x + [a \circ P_{\mathfrak{h}}, P_{\mathfrak{h}}] = \left. \Omega(a \circ P^+_{\mathfrak{m}}, P^+_{\mathfrak{m}})\right|_e.
\end{equation}

Equation \eqref{eq:bundle_vertical_curvature} demonstrates that the dynamic coupling between the Poisson leaf and the transversal subspace is governed directly by the curvature of the principal $H$-bundle. By the Frobenius theorem, the curvature $2$-form $\Omega$ measures the obstruction to the complete integrability of the horizontal distribution $\mathcal{H}$. This yields a rigid geometric reduction: when the canonical connection is flat ($\Omega \equiv 0$)—which occurs algebraically if and only if the transversal subspace $\mathfrak{m}$ forms a Lie subalgebra, $[\mathfrak{m}, \mathfrak{m}] \subseteq \mathfrak{m}$—the horizontal distribution $\mathcal{H}$ is involutive, the Lie group $G$ splits locally as a direct product $H \times M$ (where $M$ is the integral submanifold integrating $\mathfrak{m}$), and the fiber dynamics along $H$ decouple completely from the transversal evolution.
\end{rem}

In the classical theory of partial differential equations, the complete decoupling of a system into autonomous components is often treated as a purely analytical feature. However, equation \eqref{eq:bundle_vertical_curvature} establishes that for hydrodynamic flows on F$_{\text{man}}$-Lie groups, this decoupling is governed entirely by a rigid geometric obstruction: the non-linear PDE separates into independent flows if and only if the principal $H$-bundle admits a flat canonical connection.\\

To illustrate the preceding results, we present the following example.
\begin{ejem}
Let $\mathfrak{g} = \operatorname{span}\{a, b, c, d\} = \mathbb{R}a \ltimes \mathfrak{n}$, where $\mathfrak{n}$ denotes the three-dimensional Heisenberg Lie algebra, endowed with the solvable Lie algebra structure (see \cite[Theorem 1.5]{ABD}) given by the non-zero brackets:
\begin{equation*}
[a, c] = c, \qquad [a, d] = -d, \qquad [c, d] = b,
\end{equation*}
which implies that $b$ is strictly central, i.e., $[b, \cdot] = 0$.\\

We equip $\mathfrak{g}$ with a commutative and associative product $\circ$ defined by the non-zero relations:
\begin{equation*}
b \circ b = b, \qquad b \circ c = c, \qquad b \circ d = d.
\end{equation*}
The triple $(\mathfrak{g}, [\cdot, \cdot], \circ)$ forms an F$_{\text{man}}$-algebra. Moreover, this structure is non-Poisson ($L \neq 0$), as evidenced by the following non-vanishing evaluations of the algebraic Leibnizator:
\begin{equation*}
L(c, a, b) = c, \qquad L(d, a, b) = -d, \qquad L(c, c, d) = -c, \qquad L(d, c, d) = d.
\end{equation*}

We define the direct sum decomposition $\mathfrak{g} = \mathfrak{h} \oplus \mathfrak{m}$ by setting
\begin{equation*}
\mathfrak{h} = \operatorname{span}\{a, b \}, \qquad \mathfrak{m} = \operatorname{span}\{c, d \}.
\end{equation*}
This decomposition satisfies:
\begin{enumerate}[(i)]
\item For any $a \in \mathfrak{h}$ and $b, c \in \mathfrak{g}$, we have $L(a, b, c) = 0$; in particular, $L\vert{}_{\mathfrak{h}} \equiv 0$. Because $[\mathfrak{h}, \mathfrak{h}] = \{0\}$ and $\mathfrak{h} \circ \mathfrak{h} = \mathbb{R}b \subseteq \mathfrak{h}$, the subspace $\mathfrak{h}$ is a Poisson subalgebra.
\item Since $[a, c] = c \in \mathfrak{m}$ and $[a, d] = -d \in \mathfrak{m}$, it follows that $[\mathfrak{h}, \mathfrak{m}] = \mathfrak{m}$.
\item The multiplication preserves the transversal subspace, as $\mathfrak{h} \circ \mathfrak{m} = \mathbb{R}c \oplus \mathbb{R}d = \mathfrak{m}$.
\end{enumerate}Furthermore, this splitting satisfies the symmetric space condition $[\mathfrak{m}, \mathfrak{m}] = \mathbb{R}b \subseteq \mathfrak{h}$. By Remark \ref{conform}, the transversal hydrodynamic equation over $\mathfrak{m}$ linearizes completely. However, the coupling between the base flow and the Poisson leaf is governed by the curvature $2$-form, which is non-vanishing: $\left. \Omega(c, d)\right|_e = -[c, d]_{\mathfrak{h}} = -b \neq 0$.
\end{ejem}

%\vspace*{0.5cm}{\color{purple} \Large ATG: Voy aqu\'i}\\

%\section{Acknowledgments}

%%%%%%%%%%%%%%%%%%%%%%%%%%%%%%%%%%% Section

\end{document}